\documentclass{amsart}

\usepackage{amsmath,amsthm,amssymb,amscd,enumerate,mathtools,enumitem,tikz-cd,bm}
\usepackage{amsfonts}
\usepackage{rotating}
\usepackage{euscript}
\usepackage{epsfig,verbatim}

\def\be{\begin{equation}}
\def\ee{\end{equation}}

\def\C{{\mathbb C}} 
\def\f{\mathcal}
 
\def\P{{\mathbb P}}

\def\phi{{\varphi}}

\def\deg{{\rm deg\,}}

\def\GCD{{\rm gcd }}

\def\Aut{{\rm Aut}}

\def\bp{\begin{proposition}}
\def\ep{\end{proposition}}

\def\bt{\begin{theorem}}
\def\et{\end{theorem}}
\def\br{\begin{remark}}
\def\er{\end{remark}}
\def\be{\begin{equation}}
\def\bee{\begin{equation*}}
\def\l{\label}
\def\la{\label}

\def\ee{\end{equation}}
\def\eee{\end{equation*}}
\def\bl{\begin{lemma}}
\def\el{\end{lemma}}
\def\bc{\begin{corollary}}
\def\ec{\end{corollary}}
\def\pr{\noindent{\it Proof. }}

\def\bd{\begin{definition}}
\def\ed{\end{definition}}
\def\t{\widetilde}
\def\h{\widehat}

\def\t{\widetilde }
\def\h{\widehat}

\newtheorem{theorem}{Theorem}[section]
\newtheorem{lemma}[theorem]{Lemma}
\newtheorem{definition}[theorem]{Definition}
\newtheorem{corollary}[theorem]{Corollary}
\newtheorem{proposition}[theorem]{Proposition}
\newtheorem{problem}[theorem]{Problem}

\theoremstyle{definition}

\theoremstyle{definition}
\newtheorem{remark}[theorem]{Remark}
\def\bpr{\begin{problem}}
\def\epr{\end{problem}}
\def\Mon{{\rm Mon}}

\mathtoolsset{showonlyrefs}

\begin{document}

\title[]{Periodic curves for general endomorphisms\linebreak  of $\C\P^1\times \C\P^1$
}

\author[F. Pakovich]{Fedor Pakovich}
\thanks{
This research was supported by ISF Grant No. 1092/22.} 
\address{Department of Mathematics, 
Ben-Gurion University of the Negev, P.O.B. 653, Beer Sheva,  8410501, Israel}
\email{pakovich@math.bgu.ac.il}

\begin{abstract}
We show that for a general rational function $A$ of degree $m\ge 2$,
every decomposition of an iterate $A^{\circ n}$, $n\ge 1$, into
a composition of indecomposable rational functions is equivalent
to the decomposition $A\circ\cdots\circ A$ with $n$ factors.
As an application, we prove that for every $m\ge 2$ there exists
a non-empty Zariski open subset $U\subset\operatorname{Rat}_m$
such that, for any $A_1,A_2\in U$, the endomorphism of
$\mathbb{CP}^1\times\mathbb{CP}^1$ given by
$(z_1,z_2)\mapsto(A_1(z_1),A_2(z_2))$
admits an irreducible periodic algebraic curve that is neither
a vertical nor a horizontal line if and only if $A_1$ and $A_2$
are conjugate.
\end{abstract}

\maketitle

\section{Introduction}
Let $A$ be a rational function over $\C$ of degree $m \geq 2$.  
The function $A$ is said to be \emph{indecomposable} if, whenever it can be written as a composition  
$A = A_2 \circ A_1$ of rational functions, at least one of $A_1$ and $A_2$ has degree one.  
Any expression of $A$ as a composition  
\[
A = A_r \circ A_{r-1} \circ \cdots \circ A_1,
\]  
where each $A_i$ is a rational function of degree at least two, is called a \emph{decomposition} of $A$.  
Two decompositions,  
\[
A = A_r \circ A_{r-1} \circ \cdots \circ A_1 \quad \text{and} \quad A = \h A_\ell \circ \h A_{\ell-1} \circ \cdots \circ \h A_1,
\]  
are said to be \emph{equivalent} if $\ell = r$ and either $r = 1$ with $A_1 = \h A_1$, or $r \geq 2$ and there exist  
Möbius transformations $\mu_i$ for $1 \leq i \leq r - 1$ such that  
\[
A_r = \h A_r \circ \mu_{r-1}, \quad
A_i = \mu_i^{-1} \circ \h A_i \circ \mu_{i-1} \quad \text{for } 1 < i < r, \quad \text{and} \quad
A_1 = \mu_1^{-1} \circ \h A_1.
\]

In this paper, we are interested in decompositions of all iterates of a rational function $A$,  
with emphasis on the case when, for every $n \geq 1$, every decomposition of $A^{\circ n}$ into a composition of indecomposable rational functions  
is equivalent to the  decomposition $A\circ\cdots\circ A$ with $n$ factors. In this case, we say that the iterates of $A$ admit no non-trivial decompositions.  
Note that this condition implies that $A$ is itself indecomposable.

It was shown in the recent paper \cite{ite} that the iterates of a \emph{general} rational function $A$ of degree $m \geq 4$ admit no non-trivial decompositions.  
Here and below, we say that a statement holds for general rational functions of degree $m$ if,  
upon identifying the set of rational functions of degree $m$ with the algebraic variety $\mathrm{Rat}_m$  
obtained from $\mathbb{CP}^{2m+1}$ by removing the resultant hypersurface, the statement holds for all $F \in \operatorname{Rat}_m$, with the exception of some proper Zariski closed subset. 
In more detail, it was shown in \cite{ite} that the iterates of a rational function $A$ of degree $m \geq 4$ admit no non-trivial decompositions whenever $A$ is \emph{simple},  
that is, whenever for every $z \in \mathbb{CP}^1$, the preimage $A^{-1}\{z\}$ contains at least $m - 1$ points.

Although it was shown in \cite{ite} that the iterates of certain simple rational functions of degrees two and three can admit non-trivial decompositions, the main result of this paper shows that—by strengthening the simplicity condition—one can ensure that the iterates of general rational functions of degrees two and three still admit no non-trivial decompositions. This extends the result of \cite{ite}, as stated in the following theorem.

\bt \l{t1} For every \( m \geq 2 \), the iterates of a general rational function $A$ of degree $m$ admit no non-trivial decompositions.
\et 
As in \cite{ite}, Theorem~\ref{t1} can be applied to the problem of describing periodic algebraic curves for endomorphisms of \( (\mathbb{CP}^1)^2 \) of the form  
\[
(A_1, A_2) : (z_1, z_2) \mapsto (A_1(z_1), A_2(z_2)),
\]  
where \( A_1 \) and \( A_2 \) are rational functions.  
Specifically, the classification of periodic curves obtained in \cite{ic}, which incorporates earlier results for the polynomial case from \cite{ms}, relates this problem to the study of functional decompositions of iterates of \( A_1 \) and \( A_2 \).  
This connection allows Theorem~\ref{t1} to be applied to the analysis of periodic curves for general rational functions \( A_1 \) and \( A_2 \), and leads to the following result, which was proved in \cite{ite} under the more restrictive assumption \( m \geq 4 \).

\begin{theorem} \l{cur} 
For every $m \geq 2$, there exists a non-empty  Zariski open subset $U \subset \mathrm{Rat}_m$ such that the following holds. 
For any $A_1, A_2 \in U$, an irreducible algebraic curve $C \subset (\C\P^1)^2$ that is not a vertical 
or horizontal line is $(A_1, A_2)$-periodic if and only if
\[
A_2 = \alpha \circ A_1 \circ \alpha^{-1}
\]
for some Möbius transformation $\alpha$, and $C$ is one of the graphs
\[
y = (\alpha \circ A_1^{\circ s})(x) \quad \text{or} \quad x = (A_1^{\circ s} \circ \alpha^{-1})(y),
\]
for some $s \geq 0$. In particular, every such  $(A_1, A_2)$-periodic curve is $(A_1, A_2)$-invariant.
\end{theorem}

Theorem~\ref{cur} has an interesting application in complex dynamics. Namely, Zhuchao Ji and Junyi Xie recently proved in \cite{xie} that a general rational function \( A \) of degree \( m \geq 2 \) is uniquely determined up to conjugacy by its multiplier spectrum (see also \cite{mm} for an exposition and a simplified proof,
and \cite{hug} for a stronger result in the polynomial case). One of the steps in their proof relies on Theorem~\ref{cur}. Since Theorem~\ref{cur} was established in \cite{ite} only for \( m \geq 4 \), the argument in \cite{xie} does not apply to the cases \( m =2 \) and \( m =3 \). While for these degrees the result of \cite{xie} follows
from earlier work---\cite{mil1} for \(m=2\) and
\cite[Theorem~1.2]{go} for \(m=3\)
(see also \cite{hut,we})---it remained an interesting question whether the argument in \cite{xie} could be extended to \( m =2 \) and \( m =3 \) by strengthening the result of \cite{ite} to the form stated in Theorem~\ref{cur}. Addressing this question was one of the main motivations for writing the present paper.

Taking into account the results of~\cite{ite}, the main contribution of this paper is the proof of Theorems~\ref{t1} and~\ref{cur} for the remaining cases \( m =2 \) and \( m =3 \). Below, we briefly describe the explicit conditions on \( A \) under which Theorem~\ref{t1} holds for these values of $m$, along with our approach to the proof, which differs in these two cases.

Let \( A \) be a rational function of degree \( m \geq 2 \). We define \( \Sigma(A) \) as the group of Möbius transformations \( \mu \) satisfying \( A \circ \mu = A \), and the group \( \Sigma_\infty(A) \) as the union
\[
\Sigma_\infty(A) = \bigcup_{k=1}^\infty \Sigma(A^{\circ k}).
\]
Unless $A$ is conjugate to $z^{\pm m}$, the group
$\Sigma_\infty(A)$ is finite and, in some cases, can be computed explicitly (see \cite{sym}). Moreover, Theorem 2.4 in \cite{ite} implies that for a general rational function \( A \) of degree \( m \geq 3 \), the group \( \Sigma_\infty(A) \)  is trivial. However, this is not the case for \( m=2 \), since for any rational function \( A \) of degree two, the group \( \Sigma(A) \) is a cyclic group of order two, generated by some Möbius transformation which we denote by \( \mu_A \).

 Our proof of Theorem~\ref{t1} for $m=2$ proceeds as follows. First, we show  that if $A$ is a rational function of degree two, then the equality 
\be \l{eqa} \Sigma_\infty(A)=\Sigma(A)\ee
 holds whenever 
\be \l{emp}
\mu_A(V(A)) \cap V(A) = \emptyset,
\ee
where, here and below, $V(A)$ denotes the set of critical values of a rational function  $A$  (Theorem~\ref{t41}). The next step relies on the following result, which is of independent interest.

\bt \l{tt1}
Let $A$ be a rational function of degree $m \geq 2$. Then for any decomposition
\be \l{de}
A^{\circ n} = A_r \circ A_{r-1} \circ \cdots \circ A_1
\ee
of an iterate $A^{\circ n}$, $n \geq 1$, into a composition of indecomposable rational functions, the inequality $\deg A_1 \leq m$ holds.
\et

Applying this theorem to \eqref{de} when $\deg A=2$,
we see that $\deg A_1=2$ and that the group
$\Sigma_\infty(A)$ contains $\mu_{A_1}$. Under condition~\eqref{emp}, it follows that \(  \mu_{A_1}=\mu_A  \), which readily implies the conclusion of Theorem~\ref{t1} for \( A \). Finally, it is easy to show that~\eqref{emp} holds for general \( A \) (note that all functions of degree two are simple).

Now let $A$ be a simple rational function of degree three. Then it is easy to see that there exist two orbifolds $\f O_1$ and $\f O_2$ on $\mathbb{CP}^1$ with signature $\{2,2,2,2\}$ such that $A : \f O_1 \to \f O_2$ is a covering map between orbifolds. We show that the conclusion of Theorem~\ref{t1} holds whenever $\f O_1 \ne \f O_2$, that is, whenever $A$ is not a Lattès map (Theorem~\ref{l52}). Our approach here is similar to that of \cite{ite} and is based on studying the conditions under which the algebraic curve defined by the equation $$A(x) - A_r(y) = 0,$$ arising from equality \eqref{de}, may have an irreducible component of genus zero. This analysis, however, is rather complicated and involves a variety of techniques, with the main results stated in Theorems~\ref{l3}, \ref{l2}, and \ref{l51}.

This paper is organized as follows. In Section~2, we prove
Theorem~\ref{tt1}. Sections~3 and~4 are devoted to the proof
of Theorem~\ref{t1} and related results for $m=2$ and $m=3$,
respectively, following the approach outlined above.
Finally, in Section~5, we prove Theorem~\ref{cur} for these
degrees. We also establish the triviality of certain groups
of M\"obius transformations associated with general rational
functions and prove that, for a general rational function $A$
of degree $m\geq2$, the equality
\[
A^{\circ n}=B^{\circ n}
\]
for some rational function $B$ and integer $n\geq1$ implies
that $A=B$ (Theorem~\ref{it}). In particular, for each fixed
$n\geq2$, the $n$-fold iteration operator $A\mapsto A^{\circ n}$
is injective on a non-empty Zariski open subset of
$\operatorname{Rat}_m$, a result previously established by Ye
in Theorem~1.4 of~\cite{ye}. Our result provides a single
non-empty Zariski open subset for which this conclusion
holds simultaneously for all $n\geq1$.

\section{Preliminary results on decompositions of iterates}
Let $A$ and $B$ be non-constant rational functions, and let $A_1, A_2$ and $B_1, B_2$ be pairs of polynomials without common roots such that 
$$A = \dfrac{A_1}{A_2}\quad{\rm and}\quad B = \dfrac{B_1}{B_2}.$$ We define the algebraic curve $h_{A,B}$ by
\be \l{ins}
h_{A,B}:\quad A_1(x)B_2(y) - A_2(x)B_1(y) = 0.
\ee

We begin by recalling the following statement, which follows easily from general properties of fiber products (see, e.g., Section 2.1 of \cite{low}).

\bl \l{zai} Let \( A, B, X, Y \) be non-constant rational functions satisfying  
\[
A \circ X = B \circ Y \quad \text{and} \quad \mathbb{C}(X, Y) = \mathbb{C}(z).
\]  
Then  
\[
\deg Y \le \deg A, \quad \deg X \le \deg B,
\]  
and the equalities  
\[
\deg Y = \deg A, \quad \deg X = \deg B
\]  
hold if and only if the algebraic curve \( h_{A,B} \) is irreducible.
 \qed 
\el 

Theorem \ref{tt1} is obtained from Lemma \ref{zai} as follows. 

\vskip 0.2cm 
\noindent{\it Proof of Theorem \ref{tt1}.} The proof proceeds by induction on $n$. For $n = 1$, the statement is obvious. For $n > 1$, the inductive step is established by considering two cases depending on whether
\be \l{dh} 
\mathbb{C}(A^{\circ (n-1)}, A_1) = \mathbb{C}(z)
\ee
or not.

In the first case, considering the commuting diagram
\[
\begin{CD}  
\mathbb{CP}^1 @>A_1>> \mathbb{CP}^1 \\  
@VV A^{\circ (n-1)} V @VV A_r \circ A_{r-1} \circ \cdots \circ A_2 V \\  
\mathbb{CP}^1 @>A>> \mathbb{CP}^1,   
\end{CD}  
\]
we see that \( \deg A_1 \leq m \) by Lemma \ref{zai}.

In the second case, since $A_1$ is indecomposable,
L\"uroth's theorem implies that there exists a rational
function $U$ such that
\[
A^{\circ(n-1)}=U\circ A_1.
\]
If $\deg U=1$, then $A^{\circ(n-1)}$ is indecomposable,
so $n=2$ and $\deg A_1=m$.

Otherwise, write
\[
U=U_l\circ U_{l-1}\circ\cdots\circ U_1
\]
as a composition of indecomposable rational functions.
Then
\[
A^{\circ(n-1)}
=U_l\circ U_{l-1}\circ\cdots\circ U_1\circ A_1
\]
is a decomposition into indecomposable rational functions,
and the induction hypothesis gives $\deg A_1\leq m$.
\qed

Theorem \ref{tt1} implies the following corollary. 

\bc \l{c31}  
Let \( A \) be a rational function of prime degree \( p \). Then, for any decomposition 
\be \l{plk} 
A^{\circ n} =A_r \circ A_{r-1} \circ \cdots \circ A_1 
\ee 
of an iterate \( A^{\circ n} \), \( n \geq 1 \), into a composition of indecomposable rational functions, the equality \( \deg A_1 = p \) holds.
\ec
\pr Since \( \deg A = p \) is prime, equality \eqref{plk} implies that \( \deg A_1 \) is a power of \( p \). Thus, the statement follows from Theorem \ref{tt1}.
 \qed

\section{Iterates of quadratic rational functions}

\subsection{The genera of the curves $h_{A,B}$ and $h_A$}

Let $A$ and $B$ be non-constant rational functions of degrees $m$ and $l$, respectively. 
Under the assumption that the algebraic curve $h_{A,B}$ is irreducible, its genus can be calculated explicitly in
terms of the ramification of $A$ and $B$ as follows.  
Let $S = \{z_1, z_2, \ldots, z_r\}$ be the union  
of $V(A)$ and $V(B)$. For $i$, $1 \leq i \leq r$, we denote by  
$(a_{i,1}, a_{i,2}, \ldots, a_{i,p_i})$  
the collection of multiplicities of $A$ at the points of $A^{-1}\{z_i\}$, and by  
$(b_{i,1}, b_{i,2}, \ldots, b_{i,q_i})$  
the collection of multiplicities of $B$ at the points of $B^{-1}\{z_i\}$. In this notation, the following statement holds (see  \cite{f3} or \cite{low}). 
\bt  \l{t21} Let $A$ and $B$ be non-constant rational functions such that the curve $h_{A,B}$ is irreducible. Then 
\be \l{rh0}
2 - 2g(h_{A,B}) = \sum_{i=1}^r \sum_{j_2 = 1}^{q_i} \sum_{j_1 = 1}^{p_i} \gcd(a_{i,j_1}, b_{i,j_2}) - lm(r - 2). 
\ee
\et 

Below, we will use the following condition implying the irreducibility of $h_{A,B}$  in terms of $V(A)$ and $V(B)$ (see \cite{pq}, Proposition 3.1).

\bt \l{t22} Let $A$ and $B$ be non-constant rational functions such that the set 
$V(A)\cap V(B)$ contains at most one element. Then the curve $h_{A,B}$ is irreducible. \qed 
\et 

The above theorems imply the following corollary. 

\bc \l{c21} Let $A$ be a rational function of degree two,  and $\mu$ a M\"obius transformation  such that \be \l{qw} \mu(V(A)) \cap V(A) = \emptyset.\ee Then the algebraic curve $h_{A,\mu\circ A}$ is irreducible and has genus one. In particular, the functional equation 
\begin{equation} \label{ur}
A \circ X = (\mu \circ A) \circ Y
\end{equation}
has no solutions in non-constant rational functions $X$ and $Y$.

\ec 
\pr Since $$V(\mu\circ A)=\mu(V(A)),$$ it follows  from \eqref{qw} by Theorem \ref{t22} that  
the algebraic curve $h_{A,\mu\circ A}$ is irreducible. Moreover, it follows  from \eqref{qw} by Theorem \ref{t21} that  
$$2-2g(h_{A,\mu\circ A})=(1+1+1+1+1+1+1+1)-2\cdot 2\cdot 2=0.$$
Hence, $g(h_{A,\mu\circ A})=1.$ \qed 

In the case $A=B$ with $\deg A\ge2$, the algebraic
curve $h_{A,B}$ is always reducible, as it   contains the component
$x-y=0$. In this situation, it is convenient to replace the curve $h_{A,A}$ defined by \eqref{ins} with the curve  
\be \label{cur2}
h_A: \ \frac{A_1(x)A_2(y) - A_2(x)A_1(y)}{x - y} = 0.
\ee
We will use the following result (see \cite{ite}, Theorem 2.4).

\bt \l{goo1} 
Let $A$ be a simple rational function of degree $m\geq 3$. Then the curve $h_A$ 
is irreducible and  
$g(h_{A})>0$. In particular, the equality $A\circ X=A\circ Y$, where $X$ and $Y$ are non-constant rational functions, implies that $X=Y$. \qed 
\et

\subsection{\l{sub} Decompositions of iterates of quadratic rational functions}

We start by proving the following two lemmas. 

\bl \label{l41} 
Let
\[
A(z) = \frac{az^2 + bz + c}{dz^2 + ez + f}
\]
be a rational function of degree  two. Then the following statements hold:
\begin{itemize}
\item[i)] The group $\Sigma(A)$ is a cyclic group of order two generated by the Möbius transformation
\be \l{mua} 
\mu_A(z) = \frac{(cd-af)z + (ce - bf)}{(ae - bd)z + (af - cd)}.
\ee

\item[ii)] A pair of non-constant rational functions $X, Y$ satisfies the functional equation
\[
A \circ X = A \circ Y
\]
if and only if either $Y=X$ or $Y = \mu_A \circ X$.
\end{itemize}
\el
\pr The first part is obtained by directly solving $A \circ \mu = A$ in terms of the coefficients of $\mu$. To prove the second part, it is sufficient to observe that the curve $h_{A,A}$ has exactly two irreducible components: $x - y = 0$ and $\mu_A(x) - y = 0$.
 \qed

\bl \label{l42}
Let $A_1$ and $A_2$ be rational functions of degree two. Then
$\mu_{A_1} = \mu_{A_2}$ if and only if $A_2 = \nu \circ A_1$ for some M\"obius transformation $\nu$.
\el

\pr The “if” direction is immediate. To prove the “only if” direction, we proceed as follows. 
Since $$\deg A_1 = \deg A_2 = 2,$$ it follows from L\"uroth’s theorem that
$
A_2 = \nu \circ A_1
$ 
for some M\"obius transformation $\nu$ if and only if $\C(A_1, A_2) \ne \C(z)$. Hence, it suffices to show that, under the condition $\C(A_1, A_2) = \C(z)$, the identities
\be \l{st}
A_1 \circ \mu = A_1, \quad A_2 \circ \mu = A_2,
\ee
for a M\"obius transformation $\mu$ imply $\mu = \mathrm{id}$.

Assume that $\C(A_1, A_2) = \C(z)$. Then there exist polynomials $U, V \in \mathbb{C}[x, y]$ such that
\[
z = \frac{U(A_1, A_2)}{V(A_1, A_2)}.
\]
Now  substituting $\mu$ for $z$ and using \eqref{st}, we obtain
$$ 
\mu = \frac{U(A_1 \circ \mu, A_2 \circ \mu)}{V(A_1 \circ \mu, A_2 \circ \mu)}
= \frac{U(A_1, A_2)}{V(A_1, A_2)} = z. \eqno{\Box} 
$$

\bt \label{t41} 
Let $A$ be a rational function of degree two such that  
\be \l{enot} 
\mu_A(V(A)) \cap V(A) = \emptyset.
\ee
Then the group \( \Sigma_\infty(A) \)  coincides with $\Sigma(A)$.
\et 

\pr 
We prove by induction on $n$ that  
\[
\Sigma(A^{\circ n}) = \Sigma(A).
\]
For $n = 1$, the statement is trivial. Suppose now that $\nu \in \Sigma(A^{\circ n})$ for some $n > 1$. Then  
\be \label{th}
 A \circ (A^{\circ (n-1)} \circ \nu) =A \circ A^{\circ (n-1)},
\ee  
which, by Lemma \ref{l41}, implies that either  
\be \label{tha} 
A^{\circ (n-1)} = \mu_A \circ A^{\circ (n-1)} \circ \nu,
\ee 
or  
\be \label{th1}
A^{\circ (n-1)} = A^{\circ (n-1)} \circ \nu.
\ee  

In the first case, however, equation \eqref{ur} with $\mu = \mu_A$ admits a rational solution:  
\[
X = A^{\circ (n-2)}, \quad Y = A^{\circ (n-2)} \circ \nu,
\]
which implies by Corollary \ref{c21} that  
\be \label{eq}
\mu_A(V(A)) \cap V(A) \neq \emptyset,
\ee
contrary to the assumption. Thus, equality \eqref{th1} holds, and the statement follows by the induction hypothesis.  
\qed

Let us recall that writing a rational function $A = A(z)$ of degree $m$ as $A = P/Q$, where
\be \l{coe}
P(z) = a_m z^m + a_{m-1} z^{m-1} + \dots + a_0, \qquad
Q(z) = b_m z^m + b_{m-1} z^{m-1} + \dots + b_0
\ee
are polynomials without common roots, one can identify the space of rational functions of degree $m$ with the algebraic variety
\[
\mathrm{Rat}_m = \mathbb{CP}^{2m+1} \setminus \left\{ \operatorname{Res}_{m,m,z}(P,Q) = 0 \right\},
\]
where $\operatorname{Res}_{m,m,z}(P,Q)$ denotes the resultant of $P$ and $Q$. Note that 
since $\operatorname{Rat}_m$ is a non-empty Zariski
open subset of $\mathbb{CP}^{2m+1}$, it is irreducible.
Consequently, any finite intersection of non-empty
Zariski open subsets of $\operatorname{Rat}_m$
is non-empty.

For any $A \in \mathrm{Rat}_m$, the set of finite critical points of $A$ coincides with the set of zeros of its Wronskian
\be \l{w}
W(z) = P'(z)Q(z) - P(z)Q'(z).
\ee
Clearly, $\deg W \le 2m-2$, and equality holds unless $A$ lies on the projective hypersurface $L \subset \mathbb{CP}^{2m+1}$ defined by
\be \l{u}
L : a_m b_{m-1} - b_m a_{m-1} = 0.
\ee

By a standard property of the resultant, if $R(t)$ is the polynomial defined by
\be \l{res}
R(t) = \operatorname{Res}_{2m-2,\, m,\, z}\bigl(W(z),\; P(z) - Q(z)t\bigr),
\ee
then for any $A \in \mathrm{Rat}_m \setminus L$ we have
\[
R(t) = c \prod_{\zeta,\, W(\zeta)=0} \bigl(P(\zeta) - Q(\zeta)t\bigr)
\]
for some constant $c \in \mathbb{C}^*$. Hence, the zeros of $R(t)$ coincide with the finite critical values of $A$. Defining now $D$ as the projective hypersurface in $\mathbb{CP}^{2m+1}$ where the discriminant $\operatorname{Discr}(Q)$ vanishes, we see that for any $A \in \mathrm{Rat}_m \setminus (L \cup D)$ the critical points and values of $A$ are finite and $R(t)$ has degree $2m-2$.

\bl \l{sus} For a general rational function of degree two, the condition 
\eqref{enot} holds. 
\el 
\pr
Assuming that  
\[
A(z) = \frac{az^2 + bz + c}{dz^2 + ez + f} \in \operatorname{Rat}_2,
\]
and keeping the notation introduced above, we see that for $A \in \mathrm{Rat}_2 \setminus (L\cup D)$ the set of zeros of the polynomial $S(t)$ defined by  
\[
S(t) = \operatorname{Res}_{2,1,z}\bigl( R(z),\; (cd-af)z + (ce - bf) - \bigl((ae - bd)z + (af - cd)\bigr)t \bigr)
\]
coincides with the set of finite values of $\mu_A$ at the critical values of $A$.

Therefore, if $Z \subset \mathbb{CP}^{5}$ is the projective hypersurface defined by  
\[
Z : \operatorname{Res}_{2,2,t}\bigl(R(t), S(t)\bigr) = 0,
\]
then for every rational function $A$ that belongs to the set  $\operatorname{Rat}_2\setminus (L \cup D \cup Z)$ its critical values are finite  and \eqref{enot} holds. Moreover, this Zariski open subset is non-empty: for example,
it is easy to see that it contains $A(z)=z+1/z$, for which $\mu_A(z)=1/z$. \qed

\bc \l{c41} For a general rational function $A$ of degree two, the group \( \Sigma_\infty(A) \) coincides with the group \( \Sigma(A) \). 
\ec 
\pr This follows from Theorem~\ref{t41} and Lemma \ref{sus}.  
\qed

\bt \label{t42} 
Let $A$ be a rational function of degree two such that  
\[
\mu_A(V(A)) \cap V(A) = \emptyset.
\]
Then the iterates of $A$ admit no non-trivial decompositions. 
\et 

\pr 
The proof is by induction on \( n \). For \( n = 1 \), the statement holds since \( \deg A = 2 \) is prime.  

Assume now that the statement holds for $A^{\circ k}$
whenever $1\le k<n$, and let  
\be \label{e} 
A^{\circ n} = A_r \circ A_{r-1} \circ \cdots \circ A_1
\ee  
be a decomposition of \( A^{\circ n} \) into a composition of indecomposable rational functions. It is easy to see that to complete the inductive step, it suffices to show that  
\be \label{mo} 
A_1 = \mu \circ A
\ee  
for some Möbius transformation \( \mu \).  

By Corollary \ref{c31}, the function \( A_1 \) has degree two, and equality \eqref{e} implies that  
$
\mu_{A_1} \in \Sigma_\infty(A).
$  
Thus, by Theorem \ref{t41}, we must have \( \mu_{A_1} = \mu_A \), and  \eqref{mo} follows from Lemma \ref{l42}. 
\qed

\bc \l{c42}  For a general rational function $A$ of degree two, the 
iterates of $A$ admit no non-trivial decompositions. 
\ec 
\pr This follows from Theorem~\ref{t42} and Lemma \ref{sus}. \qed

For a rational function \( F \), we define the group \( G(F) \) as the group of all M\"obius transformations \( \nu \) such that  
\be \l{re} 
F \circ \nu = \delta \circ F 
\ee  
for some M\"obius transformation \( \delta \).  
Notice that  equality \eqref{re} readily implies  the equalities 
$$ 
\delta(V(F)) = V(F), \quad \nu(C(F)) = C(F), 
$$  
where \( C(F) \) denotes the set of critical points of \( F \).  

We define the group \( \operatorname{Aut}(A) \) as the subgroup of \( G(A) \) consisting of all M\"obius transformations \( \sigma \) such that  
$$ 
A \circ \sigma = \sigma \circ A, 
$$  
and define the group \( \operatorname{Aut}_\infty(A) \) by  
\[
\operatorname{Aut}_\infty(A) = \bigcup_{k=1}^\infty \operatorname{Aut}(A^{\circ k}).
\]

\bc \l{c52}  For a general rational function $A$ of degree two, the group 
$\operatorname{Aut}_\infty(A)$ is trivial. 
\ec 
\pr Let $A$ be a rational function of degree two satisfying
\eqref{enot}. We start by showing that
$\operatorname{Aut}_\infty(A)$ is isomorphic to a subgroup
of the Klein four-group $V_4$.

Let $\mu\in\operatorname{Aut}(A^{\circ k})$ for some $k>1$.
Applying Theorem~\ref{t42} to the equality
\be \l{imp}
A^{\circ k}
=\mu^{-1}\circ A^{\circ k}\circ\mu
=(\mu^{-1}\circ A)\circ A^{\circ(k-2)}\circ(A\circ\mu),
\ee
we see that there exist M\"obius transformations
$\nu$ and $\delta$ such that
\be \l{allo}
\mu^{-1}\circ A=A\circ\nu
\ee
and
\be \l{alo}
A\circ\mu=\delta\circ A.
\ee
Thus, $\mu\in G(A)$, and
\be \l{condi}
\mu(V(A))=V(A),\qquad \mu(C(A))=C(A).
\ee
For $k=1$, equality~\eqref{alo} holds with $\delta=\mu$,
and equalities~\eqref{condi} hold as well.

Since the critical points of $A$ are precisely the fixed
points of $\mu_A$, condition~\eqref{enot} implies that
$C(A)\cap V(A)=\emptyset$.
Thus, $C(A)\cup V(A)$ consists of four distinct points,
\be \l{lab}
C(A)=\{z_1,z_2\},\qquad V(A)=\{z_3,z_4\}.
\ee
By \eqref{condi}, each $\mu\in\operatorname{Aut}_\infty(A)$
induces a permutation $\sigma_\mu\in S_4$, defined by
\[
\mu(z_i)=z_{\sigma_\mu(i)},\qquad 1\le i\le4.
\]
Moreover, it is easy to see that the map $\mu\mapsto\sigma_\mu$ is a homomorphism
from $\operatorname{Aut}_\infty(A)$ into the Klein four-group
\[
V_4=\{e,\ (12)(34),\ (12),\ (34)\}\subset S_4,
\]
whose kernel is trivial, since any M\"obius transformation
fixing four distinct points is the identity.
Thus, $\operatorname{Aut}_\infty(A)$ is isomorphic
to a subgroup of $V_4$.

We now show that $\operatorname{Aut}_\infty(A)$
is not isomorphic to $V_4$, and hence
\be \l{ineq}
|\operatorname{Aut}_\infty(A)|\le2.
\ee
Assume the contrary. Then there exists a nonidentity
involution $\mu\in\operatorname{Aut}_\infty(A)$
fixing both points of $C(A)$.
Since a nonidentity M\"obius involution is uniquely
determined by its two fixed points, we have $\mu=\mu_A$. Then \eqref{enot} cannot hold, in view of the first equality
in~\eqref{condi}. Thus, \eqref{ineq} holds.

Finally, we show that inequality~\eqref{ineq} implies
the equality
\be \l{bura}
\operatorname{Aut}_\infty(A)=\operatorname{Aut}(A).
\ee
Let $\mu\in\operatorname{Aut}(A^{\circ k})$, $k>1$,
be nonidentity, and let $\delta$ satisfy \eqref{alo}.
Since $\mu$ commutes with $A^{\circ k}$, we have
\[
A^{\circ k}\circ\delta\circ A=A^{\circ k}\circ A\circ\mu
=A\circ\mu\circ A^{\circ k}
=\delta\circ A^{\circ k}\circ A.
\]
Cancelling $A$ on the right, we obtain
$\delta\in\operatorname{Aut}(A^{\circ k})$.
By \eqref{ineq}, either $\delta=\mathrm{id}$ or $\delta=\mu$.
In the first case, \eqref{alo} implies that $\mu=\mu_A$,
contradicting \eqref{condi} and \eqref{enot}.
In the second case, \eqref{alo} implies that
$\mu\in\operatorname{Aut}(A)$.
Thus, \eqref{bura} holds.

By Lemma~\ref{sus}, condition~\eqref{enot} holds
for a general rational function of degree two.
Since $\operatorname{Aut}(A)$ is trivial for a general quadratic rational function (see \cite{mil1}, Section~5), the corollary
follows. \qed

\section{Iterates of  cubic rational functions}

\subsection{Algebraic curves $h_{A,B}$ with $\deg A=3$: the irreducible case}
We begin by recalling some definitions and results concerning functional equations that employ the notion of orbifolds on the Riemann sphere.

An {\it orbifold} $\f O$ on $\C\P^1$ is a ramification function $\nu:\C\P^1\rightarrow \mathbb N$, which takes the value $\nu(z)=1$ except at a finite set of points.  
For an orbifold $\f O$, the {\it Euler characteristic} of $\f O$ is the number  
$$ \chi(\f O)=2+\sum_{z\in \C\P^1}\left(\frac{1}{\nu(z)}-1\right),$$ the set of {\it singular points} of $\f O$ is the set  
$$c(\f O)=\{z_1,z_2, \dots, z_s\}=\{z\in \C\P^1 \mid \nu(z)>1\},$$  
and the {\it signature} of $\f O$ is the multiset  
$$\nu(\f O)=\{\nu(z_1),\nu(z_2), \dots , \nu(z_s) \}.$$

We will always assume that the orbifolds under consideration are \emph{good}, meaning distinct from those with a single ramified point, or with two ramified points $z_1, z_2$ such that $\nu(z_1) \neq \nu(z_2)$ (for more details, see, e.g., \cite{lattes}).

Let $A$ be a rational function, and let $\f O_1$, $\f O_2$ be orbifolds with ramification functions $\nu_1$, $\nu_2$. We say that $A: \f O_1 \rightarrow \f O_2$ is a {\it covering map} between orbifolds if, for every $z \in \C\P^1$, the equality
\be \l{cov}
\nu_{2}(A(z)) = \nu_{1}(z) \cdot \deg_z A
\ee
holds, where $\deg_z A$ is the multiplicity of $A$ at $z$. 
If  the weaker condition  
\be \l{uuss} 
\nu_2(A(z)) \mid \nu_1(z) \deg_z A 
\ee  
holds instead of \eqref{cov} for every $z \in \C\P^1$, then we say that 
$A: \mathcal{O}_1 \to \mathcal{O}_2$  is a {\it holomorphic map} between orbifolds.

If $A: \f O_1 \rightarrow \f O_2$ is a covering map between orbifolds, then the Riemann--Hurwitz formula implies that
\be \l{rhor} \chi(\f O_1) =  \chi(\f O_2) \cdot\deg A. \ee
For holomorphic maps, the following statement holds (see \cite{semi}, Proposition 3.2).

\bp \l{p1} Let $A: \f O_1 \rightarrow \f O_2$ be a holomorphic map between orbifolds. Then
\be \l{iioopp} \chi(\f O_1) \leq \chi(\f O_2) \cdot \deg A, \ee
and equality holds if and only if $A: \f O_1 \rightarrow \f O_2$ is a covering map. \qed
\ep

Let $A$ be a non-constant rational function. If $\C\P^1$ is equipped with a ramification function $\nu_2$, then to define a ramification function $\nu_1$ on $\C\P^1$ so that $A$ becomes a holomorphic map between orbifolds $\f O_1$ and $\f O_2$,  
condition \eqref{uuss} must be satisfied, and it is easy to see that  
for any $z \in \C\P^1$, a minimal possible value for $\nu_1(z)$ is determined by the equality  
\be \l{rys} 
\nu_{2}(A(z)) = \nu_{1}(z) \GCD(\deg_z A, \nu_{2}(A(z))).  
\ee  
If \eqref{rys} holds for all $z \in \C\P^1$, we say that $A$ is a  
{\it minimal holomorphic map}  
between orbifolds $\f O_1$ and $\f O_2$. Notice that  \eqref{cov} implies \eqref{rys}. Therefore,  if $A : \f O_1 \rightarrow \f O_2$ is a covering map between orbifolds, then $A : \f O_1 \rightarrow \f O_2$ is also a minimal holomorphic map between orbifolds.

With each rational function $A$, one can naturally associate two orbifolds $\f O_1^A$ and $\f O_2^A$ by defining $\nu_2^A(z)$ as the least common multiple of the multiplicities of $A$ at the points of the preimage $A^{-1}\{z\}$, and setting  
$$
\nu_1^A(z) = \frac{\nu_2^A(A(z))}{\deg_z A}.
$$  
By construction,  
$
A:\, \f O_1^A \rightarrow \f O_2^A
$  
is a covering map between orbifolds. Orbifolds defined in this way are useful for studying the functional equation  
\be \label{mm} 
A \circ X = B \circ Y, 
\ee  
where $A$, $B$,  $X$, $Y$ are rational functions, which we usually represent by the commuting diagram  
\be \label{di1} 
\begin{CD}
\C\P^1 @>X>> \C\P^1  \\
@VV Y V @VV A V \\ 
\C\P^1  @>B>> \C\P^1. 
\end{CD}
\ee

The main result we use to analyze equation \eqref{mm} is the following statement (see \cite{semi}, Theorem 4.2 or \cite{pak}, Theorem 1.1).

 \bt \l{goodt} 
Let $A$, $B$, $X$, $Y$ be rational functions satisfying \eqref{mm} such that the curve $h_{A,B}$ is irreducible and $\C(X,Y) = \C(z)$.  
Then the  diagram  
\be \l{ddii}
\begin{CD}
\f O_1^Y @>X>> \f O_1^A \\
@VV Y V @VV A V \\ 
\f O_2^Y @>B>> \f O_2^A
\end{CD}
\ee  
consists of minimal holomorphic maps between orbifolds. \qed  
\et 

We recall that a rational function $A$ of degree $m \geq 2$ is called \emph{simple} if, for every $z \in \mathbb{CP}^1$, the preimage  
$A^{-1}\{z\}$ contains at least $m - 1$ points.  
This condition is equivalent to requiring that $A$ has exactly $2m - 2$ critical values (see \cite{ite}, Lem-\- ma 2.1).  
For a simple rational function $A$ of degree $m \geq 2$, its monodromy group satisfies 
\be \l{mon} 
\operatorname{Mon}(A) \cong S_m
\ee
(see, e.g., \cite{ite}, Theorem 2.2).

In this section, we are concerned with simple rational functions of degree three. Clearly, if $A$ is such a function, then
\be \l{nu} 
\nu(\f O_2^A) = \{2,2,2,2\}, \qquad \nu(\f O_1^A) = \{2,2,2,2\}.
\ee
In particular,
\[
\chi(\f O_2^A) = 0, \qquad \chi(\f O_1^A) = 0.
\]

\bt \l{l3} 
Let $A$, $B$, $X$, $Y$ be rational functions satisfying $A \circ X = B \circ Y$ such that the curve $h_{A,B}$ is irreducible and $\C(X,Y) = \C(z)$. Assume that $A$ is a simple rational function of degree three, $\deg B\geq2,$  and $\deg B \neq 2, 4$. Then the equalities
\[
\f O_2^A = \f O_2^B, \qquad \f O_1^A = \f O_2^X
\]
hold. 
\et  

\pr 
Let us first show that
\be \l{lu}
\nu(\f O_2^Y) = \{2, 2, 2, 2\}.
\ee
Since $B : \f O_2^Y \rightarrow \f O_2^A$ is a minimal holomorphic map between orbifolds by Theorem~\ref{goodt}, it follows from the first equality in \eqref{nu} and the definition of a minimal holomorphic map \eqref{rys} that 
\[
\nu(\f O_2^Y) = \underbrace{\{2, 2, \ldots, 2\}}_{\ell}
\]
for some natural number $\ell$. Moreover, $\ell \geq 2$, since any rational function $Y$ of degree at least two has at least two critical values, and it is easy to see that in fact $\ell = 4$. Indeed, if $2 \leq \ell \leq 3$, then $\chi(\f O_2^Y) > 0$, contradicting \eqref{iioopp} because $\chi(\f O_2^A) = 0$. On the other hand, the inequality $\ell > 4$ is impossible, since $\deg Y = \deg A = 3$ by Lemma~\ref{zai}, and a rational function of degree three has at most four critical values. Thus, \eqref{lu} holds. Notice that this implies by Proposition~\ref{p1} that 
$
B : \f O_2^Y \rightarrow \f O_2^A
$
is a covering map between orbifolds.

Now we show that
\be \l{td} c(\f O_2^B) \subseteq c(\f O_2^A).
\ee
Since $c(\f O_2^B) = V(B)$ and $c(\f O_2^A) = V(A)$, it follows easily from the Riemann--Hurwitz formula that the last inclusion is equivalent to the equality
\[
|B^{-1}(c(\f O_2^A))| = (|c(\f O_2^A)| - 2) \deg B+2 = 2 \deg B + 2,
\]
which in turn follows from the fact that $B : \f O_2^Y \rightarrow \f O_2^A$ is a covering map between orbifolds. Indeed, taking into account the equalities \eqref{nu} and \eqref{lu} for the signatures of $\f O_2^Y$ and $\f O_2^A$, it follows from \eqref{cov} that the multiplicity of $B$ at points of $B^{-1}(c(\f O_2^A))$ equals two, except at four points where it equals one. Hence, the number of points in $B^{-1}(c(\f O_2^A))$ equals 
\[ 
\frac{4 \cdot \deg B - 4}{2} + 4 = 2 \deg B + 2.  
\] 

Since $\deg B\geq 2$, to deduce the equality $\f O_2^A = \f O_2^B$ from \eqref{td}, we need only show that the equalities $\nu(\f O_2^B) = \{2,2,2\}$ or $\nu(\f O_2^B) = \{2,2\}$ are impossible. However, it is well known and easily seen that in these cases, up to replacing $B$ with $\mu_1 \circ B \circ \mu_2$ for some Möbius transformations $\mu_1$ and $\mu_2$, the function $B$ is equal to $\frac{1}{2}(z^2 + 1/z^2)$ or $z^2$, implying that $\deg B$ is either four or two. 

Finally, to prove the equality $\f O_1^A = \f O_2^X$, we observe that switching the roles of $A$ and $B$ and applying Theorem~\ref{goodt} again, we see that 
$
A : \f O_2^X \rightarrow \f O_2^B
$
is a minimal holomorphic map between orbifolds. Since $\f O_2^B = \f O_2^A$, it now follows easily from the definition of $\f O_2^A$ and formula \eqref{rys} that 
$
\f O_2^X = \f O_1^A
$. 
\qed

\subsection{Algebraic  curves $h_{A,B}$ with $\deg A=3$: the reducible case}
Let us recall that a holomorphic map between compact Riemann surfaces
 $F:\,  N\rightarrow  R$   is called a {\it Galois covering} if its group of covering transformations  
$$\Aut( N,F)=\{\sigma\in \Aut( N)\,:\, F\circ\sigma=F\}$$ acts transitively on  fibers  of $F$.  Denoting by $\f M(R)$ the field of meromorphic functions on a compact Riemann surface $R$, we can restate this condition as the condition that 
the field extension 
$\f M( N)/F^*\f M( R)$ is a Galois extension. 

In case $F$ is a Galois  covering, for the corresponding Galois group the isomorphism  
\be \l{ga} {\rm Gal}\left( \f M( N)/F^* \f M( R)\right)\cong \Aut( N,F)\ee holds.   
Notice that 
$F$ is a Galois covering 
if and only if the equality 
\be \l{dega} \vert \Aut( N,F)\vert =\deg F\ee 
holds.

Let $A$  be  a rational function.  Then the {\it normalization} of $A$ is defined as a compact Riemann surface  $N_A$ together with a holomorphic Galois covering  of the 
lowest possible degree $\h  A:N_A\rightarrow \C\P^1$ such that
 $$\h  A=A\circ H$$ for some  holomorphic map $H:\,N_A\rightarrow \C\P^1$. 
The map $\h  A$ is defined up to the change $\h  A\rightarrow \h  A\circ \alpha,$ where $\alpha\in\Aut(N_A)$, and is 
characterized by the property that  the field extension 
$\f M(N_A)/{\h  A}^*\f M(\C\P^1)$ is isomorphic to the Galois closure $\t {\f M( \C\P^1)}/A^*\f M(\C\P^1)$
of the extension $\f M( \C\P^1)/A^*\f M(\C\P^1)$.  Note that  the corresponding Galois group  is isomorphic to
the monodromy group of $A$: 
 $${\rm Gal}\big(\t{\f M(\C\P^1)}/A^*\f M(\C\P^1)\big)\cong \Mon(A)$$ 
(see, e.g., \cite{har}). In particular, this implies by \eqref{ga} and \eqref{dega} that 
\be \l{h} \vert  \Mon(A)\vert =\deg \h  A.\ee

The main technical tool for working with reducible curves $h_{A,B}$ is the following result of Fried (see \cite{f2}, Proposition 2, or \cite{pak}, Theorem 6.16).  

\bt \l{fr} Let $A$ and $B$ be non-constant rational functions such that the curve $h_{A,B}$ is reducible. Then there exist rational functions $A_1,$ $A_2$ and $B_1,B_2$ such  that
\be \la{e6} A= A_1\circ A_2, \ \ \ B=B_1\circ B_2,\ee the number of components of $h_{A_1,B_1}$ is equal to the number of components of $h_{A,B}$,  
and  $\h A_1=\h B_1$. \qed 
\et 

Note that both $A_1$ and $B_1$ must have degree at least two since otherwise the curve $h_{A_1,B_1}$ is irreducible.

\bt \l{l2} 
Let $A$ be a simple rational function of degree three, and let $B$ be a rational function such that the curve $h_{A,B}$ is reducible. Then  
\be \l{iu} 
B = A \circ R
\ee  
for some rational function $R$. In particular, $V(A) \subseteq V(B)$. 
\et  

\pr 
 Since any rational function of degree three is indecomposable, it follows from Fried's theorem that there exist rational functions $U$ and $D$ such that 
\be \l{ui} 
B = U \circ D,
\ee  
the curve $h_{A,U}$ is reducible, and $\h{A} = \h{U}$. The last equality implies that there exist holomorphic maps $E, F: N_A \to \mathbb{C}\mathbb{P}^1$ such that  
\be \l{oi} 
\h{A} = A \circ E = U \circ F.\ee Furthermore, since $A$ is simple, its monodromy group is $S_3$ by \eqref{mon}. Therefore, \(\deg \h{A} = |S_3| = 6\) by \eqref{h}, and hence \(\deg U\) equals \(6\), \(2\), or \(3\).

Since the curve \(h_{A,U}\) is always irreducible when the degrees of \(A\) and \(U\) are coprime (see, e.g., \cite{pak}, Proposition 6.7), the case \(\deg U = 2\) is impossible. Moreover, if \(\deg U = 6\), then \(\deg F = 1\) and 
\[
U = A \circ E \circ F^{-1},
\]
which together with \eqref{ui} implies \eqref{iu}.

Assume now that  $\deg U = 3$ and hence $\deg E = \deg F = 2$.
As $\widehat{A}$ is a Galois covering, it follows from \eqref{oi} that there exist subgroups $\Gamma_1$ and $\Gamma_2$ of order two in $\Aut(N_A, \widehat{A})$ such that  
\be \label{iff}
\Gamma_1 = \Aut(N_A, E), \quad \Gamma_2 = \Aut(N_A, F).
\ee  
However, all subgroups of order two in $S_3$ are conjugate, and hence there exists $\mu \in \operatorname{Aut}(N_A,\widehat A)$ such that  
\[
\Gamma_2 = \mu^{-1} \circ \Gamma_1 \circ \mu,
\]  
implying that  
\[
E \circ \mu = \delta \circ F
\]  
for some Möbius transformation $\delta$. Since $\widehat{A} = \widehat{A} \circ \mu$ for any $\mu \in \operatorname{Aut}(N_A,\widehat A)$, we obtain  
\[
\widehat{A} = A \circ E = A \circ E \circ \mu = A \circ \delta \circ F.
\]  
Since on the other hand $\widehat{A} = U \circ F$, we conclude  that  $U = A \circ \delta,$   
which, combined with \eqref{ui}, yields the required equality \eqref{iu}. 
\qed

\subsection{Decompositions of iterates of cubic rational functions} 
We recall that a \textit{Lattès map} can be defined as a rational function $A$ of degree at least two such that there exists an orbifold $\f O$ for which
$A : \f O \rightarrow \f O$
is a covering map between orbifolds (see \cite{mil2}, \cite{lattes}). Note that for such $\f O$, necessarily $\chi(\f O) = 0$ by \eqref{rhor}.

\bt \label{l51}
Let $A$ be a simple rational function of degree three. Assume that there exists a rational function $B$ such that $\deg B\geq 2,$  $\deg B \ne 2, 4$, the curve $h_{A,B}$ is irreducible, and the curve $h_{A^{\circ 2},B}$ has  an irreducible component of genus zero. Then $A$ is a Lattès map.
\et
\pr The condition that $h_{A^{\circ 2},B}$ has an irreducible component of genus zero is equivalent to the existence of rational functions $X$ and $Y$ such that $\C(X, Y) = \C(z)$ and the diagram  
\[
\begin{CD}  
\mathbb{CP}^1 @> Y >> \mathbb{CP}^1 \\  
@V X VV @V B VV \\  
\mathbb{CP}^1 @> A^{\circ 2} >> \mathbb{CP}^1  
\end{CD}  
\]
commutes.
By the universal property of fiber products, this diagram can be extended to  
\be \l{mi}  
\begin{CD}
\mathbb{CP}^1 @> Y_2 >> \mathbb{CP}^1 @> Y_1 >> \mathbb{CP}^1 \\
@V X VV @V E VV @V B VV \\
\mathbb{CP}^1 @> A >> \mathbb{CP}^1 @> A >> \mathbb{CP}^1,
\end{CD}
\ee  
where $Y_1,Y_2,$ and $E$ are rational functions satisfying $Y = Y_1 \circ Y_2$ and $$ \C(E,Y_1) = \C(z), \quad \C(X,Y_2) = \C(z).$$

Applying Theorem \ref{l3} to the right square in diagram \eqref{mi}, we see that 
\be \l{si}  
\f O_1^A= \f O_2^E.  
\ee  
Thus, since  
$A : \f O_1^A \to \f O_2^A$  
is a covering map between orbifolds, to prove the theorem it suffices to show  that 
\be \l{si2}  
\f O_2^A=\f O_2^E.  
\ee  

By Lemma~\ref{zai} applied to the right square in diagram~\eqref{mi},
we have \linebreak $\deg E=\deg B$. Hence $\deg E\geq2$ and
$\deg E\notin\{2,4\}$. Therefore, if the curve $h_{A,E}$
is irreducible, equality~\eqref{si2} follows from Theorem~\ref{l3} 
applied to the left square.

 On the other hand,  
if $h_{A,E}$ is reducible, then $V(A) \subseteq V(E)$ by Theorem \ref{l2}.  
Since \eqref{si} implies by \eqref{nu} that $$\nu(\f O_2^E)=\nu(\f O_2^A)=\{2,2,2,2\},$$  this can only happen if  
\eqref{si2} holds.  
\qed
\bt \l{l52}  
Let $A$ be a simple rational function of degree three that is not a Latt\`es map.  
Then  the iterates of $A$ admit no non-trivial decompositions.
\et  

\pr The proof is by induction on \( n \). For \( n = 1 \), the statement holds since \( \deg A = 3 \) is prime.  
Now let  
\be \l{en}  
A^{\circ n} = A_r \circ A_{r-1} \circ \cdots \circ A_1  
\ee  
be a decomposition of $A^{\circ n}$ into a composition of indecomposable rational functions, with $n \geq 2$. Since Theorem~\ref{goo1} implies that the equality $A \circ X = A \circ Y$ yields $X = Y$, to prove the inductive step it suffices to show that \eqref{en} implies $A_r = A \circ \mu$ for some Möbius transformation $\mu$. 

Clearly, \eqref{en} implies that the algebraic curve
\be \label{us}
A^{\circ 2}(x) - A_r(y) = 0
\ee
has an irreducible component of genus zero. Since $A$ is not a Lattès map and \eqref{en} implies that $\deg A_r$ is a power of $3$ (hence not $2$ or $4$), it follows from Theorem~\ref{l51} that $h_{A,A_r}$ is reducible. Therefore, the equality $A_r = A \circ R$ holds for some rational function $R$ by Theorem~\ref{l2}. Finally,  since $A_r$ is indecomposable,  $\deg R=1$. 
\qed

Since it is well known that a general rational function of degree $m \geq 2$ is not a Latt\`es map, Theorem~\ref{l52} implies the following corollary.

\bc \l{c51} For a general rational function $A$ of degree three, the iterates of \( A \) admit no non-trivial decompositions. \qed
\ec

Since Theorem~\ref{t1} is proved in \cite{ite} for $m \geq 4$, Corollary~\ref{c51}, combined with Corollary~\ref{c42}, completes the proof of Theorem~\ref{t1}.

The following result is a direct corollary of Theorem~\ref{t1}.

\begin{corollary}
\label{cor1}
For a general rational function $A$ of degree $m\geq2$,
the following holds: whenever $r\geq2$ and $G_i$,
$1\leq i\leq r$, are rational functions of degree at least
two satisfying
\[
A^{\circ n}=G_r\circ G_{r-1}\circ\cdots\circ G_1
\]
for some $n\geq1$, there exist M\"obius transformations
$\nu_i$, $1\leq i<r$, and integers $s_i\geq1$,
$1\leq i\leq r$, such that
\[
G_r=A^{\circ s_r}\circ\nu_{r-1},\quad
G_i=\nu_i^{-1}\circ A^{\circ s_i}\circ\nu_{i-1},
\quad 1<i<r,\quad\text{and}\quad
G_1=\nu_1^{-1}\circ A^{\circ s_1}.
\]
\end{corollary}

\begin{proof}
To prove the corollary, it suffices to decompose each
$G_i$, $1\leq i\leq r$, into a composition of indecomposable
rational functions and then apply Theorem~\ref{t1}.
\end{proof}

\section{Decompositions of iterates and periodic curves}
\subsection{Generalized Lattès maps} 
To study the functional equation
\[
A \circ X = X \circ B,
\]
in rational functions, as well as invariant algebraic curves for endomorphisms
\[
(A_1, A_2) : (z_1, z_2) \mapsto (A_1(z_1), A_2(z_2))
\]
of $(\mathbb{C}\mathbb{P}^1)^2$, where $A_1$ and $A_2$ are rational functions, it is useful to make use of \emph{generalized Latt\`es maps}, which generalize ordinary Latt\`es maps.

A generalized Lattès map can be defined as a rational function $A$ of degree at least two for which there exists an orbifold $\f O$, distinct from the non-ramified sphere, such that $A: \f O \to \f O$ is a minimal holomorphic map. 
Thus, 
$A$ is a Lattès map if, for some orbifold with the ramification  function $\nu$, 
\be \l{sin}
\nu(A(z)) = \nu(z) \deg_z A, \quad z \in \C\P^1,
\ee
while $A$ is a generalized Lattès map if 
\be \l{ee}
\nu(A(z)) = \nu(z) \gcd(\deg_z A, \nu(A(z))), \quad z \in \C\P^1.
\ee 
Since \eqref{sin} implies \eqref{ee}, every Lattès map is a generalized Lattès map. 

Note that \eqref{rhor} and \eqref{iioopp} imply that if $
A: \f O \to \f O
$ is a covering map, then $\chi(\f O)=0$ and 
the signature $\nu(\f O)$  belongs to the following well-known list 
\be \l{o1} 
 \{2, 2, 2, 2\},\qquad \{3, 3, 3\},\qquad \{2, 4, 4\},\qquad \{2, 3, 6\}.
\ee
On the other hand, if $
A: \f O \to \f O
$ is a minimal holomorphic  map but not a covering map, then $\chi(\f O)>0$ and
$\nu(\f O)$ belongs to the list 
\[
\begin{aligned}
\{n, n\},\; n \geq 2,\qquad \{2, 2, n\},\; n \geq 2,\qquad \{2, 3, 3\}, \qquad\{2, 3, 4\},\qquad \{2, 3, 5\}.
\end{aligned}
\]

The concept of a generalized Lattès map is useful for two  reasons. First, excluding such maps allows for simpler formulations of results concerning semiconjugate rational functions and invariant curves (see Theorems~\ref{tt} and~\ref{1} below). Second, a general rational function is not a generalized Lattès map.

In this section, we prove the latter statement. For \(m \geq 4\), it follows from Lemma~4.3 in~\cite{ite}, and below we extend it to all \(m \geq 2\). 
We begin with the following result.

\bl \l{d2} Let $A$ be a rational function of degree two, and let $\f O$ be a good  orbifold of positive Euler characteristic, distinct from the non-ramified sphere, such that $A: \f O \rightarrow \f O$ is a minimal holomorphic map. Then
\be \l{stt}
A(V(A)) \cap V(A) \neq \emptyset
\ee
unless $\nu(\f O) = \{2,2\}$.
\el 

\pr 
Let us first show that if $A: \f O \rightarrow \f O$ is a minimal holomorphic map, then 
the signature of $\f O$ cannot be $\{2, 3, 5\}$, $\{2, 3, 4\}$, $\{2, 3, 3\}$, or $\{2,2,2\}$.  

Since $A$ has exactly two critical values,
if $c(\f O)$ contains three points, then at least one
of them, say $z_0$, is not a critical value of $A$. Consequently, by \eqref{ee}, at the two distinct points $\{z_1, z_2\} = A^{-1}(z_0)$, the function $\nu$ takes the same value, namely $\nu(z_0)$. This rules out the signatures $\{2, 3, 4\}$ and $\{2, 3, 5\}$, which do not contain two equal entries.

To rule out the signature $\{2,2,2\}$, write
$c(\mathcal O)=\{z_1,z_2,z_3\}$. Observe that, by \eqref{ee},
for every point $z_0$ with $\nu(z_0)=2$, its preimage
under $A$ either consists of two points where $\nu$
takes the value $2$ (if $z_0$ is not a critical value
of $A$) or contains no such points at all (if $z_0$
is a critical value of $A$). Hence, the number of points
of $c(\mathcal O)$ in $A^{-1}\{z_1,z_2,z_3\}$ is even.
On the other hand, \eqref{ee} implies that
$\{z_1,z_2,z_3\}\subseteq A^{-1}\{z_1,z_2,z_3\}$,
so this number equals three, a contradiction. 

Now, if $\nu(\f O) = \{2, 3, 3\}$, then \eqref{ee} implies that a point $z_0$ with $\nu(z_0)=2$ satisfies $A(z_0)=z_0$ and $\deg_{z_0} A = 1$. Thus,
$z_0$ is not a critical value of $A$, and at the second point of $A^{-1}(z_0)$ the function $\nu$  also takes the value $2$, a contradiction. Hence, this signature is also impossible. 

To finish the proof, we must show that \eqref{stt} holds whenever $A: \f O \rightarrow \f O$ is a minimal holomorphic map for an orbifold $\f O$ with signature $\{2,2,n\}$ or $\{n,n\}$, where $n \geq 3$. In the first case, the point $z_0$ with $\nu(z_0)=n$ satisfies $A(z_0)=z_0$. Moreover, either $\deg_{z_0}A = 1$, or $\deg_{z_0}A = 2$ and then $\gcd(n,2)=1$. The first subcase is impossible; for otherwise, as above, we would conclude that there exist two points where $\nu$ takes the value $n$. In the second subcase, $z_0$ is a fixed point
that is also a critical value of $A$,
and hence \eqref{stt} holds.

Finally, if $\nu(\f O) = \{n,n\}$ with $n \geq 3$, and $z_1$ and $z_2$ are points such that \linebreak $\nu(z_1)=\nu(z_2)=n$, then, using arguments similar to those above, one can easily see that   $A(\{z_1,z_2\})=\{z_1,z_2\}.$ Moreover,  $n$ is odd, and the points $z_1$ and $z_2$ are critical values of $A$. Thus, \eqref{stt} holds. \qed

\bt \label{lats}  
For every \( m \geq 2 \), a general rational function \( A \) of degree $m$ is not a generalized Latt\`es map.  
\et
\pr 
Since for every $m\geq 2$ there exists a proper closed subset $N$  of $\operatorname{Rat}_m $ such that any function  $A\in \operatorname{Rat}_m\setminus N$ is simple with finite critical values  (see Lemma 3.9 in \cite{ite}), it suffices to consider only such functions. 

For a finite subset $S \subset \C\P^1$, define $A^{-1}_1\{S\}$ as the subset of $A^{-1}\{S\}$ consisting of points at which the multiplicity of $A$ equals one.
It is easy to see that if $A: \f O \rightarrow \f O$ is a minimal holomorphic map, then  \eqref{ee} implies 
\be \l{bn} 
A^{-1}_1\{c(\f O)\}\subseteq c(\f O), \qquad A\{c(\f O)\}\subseteq c(\f O).
\ee
If $A$ is a simple rational function of degree $m \geq 4$, then for any  $S$ the inequality 
\[
|A^{-1}_1\{S\}| \geq (m - 2)|S| \geq 2|S|
\]
holds, making the first inclusion in \eqref{bn} impossible. 
Thus, a simple rational function of degree $\geq 4$ cannot be a generalized Latt\`es map.

Furthermore, for a simple rational function $A$ of degree $3$ we have 
$
|A^{-1}_1\{S\}| \geq |S|,
$ 
and the equality is attained if and only if $S\subseteq V(A).$ Thus, the first inclusion in \eqref{bn} implies that \(c(\f O) \subseteq V(A)\). It then follows from the second inclusion that if a simple rational function of degree $3$ is a generalized Latt\`es map, then \eqref{stt} holds. Hence, to prove the theorem for $m = 3$, it suffices to show that for a general rational function $A$ of degree $3$,
\be \l{av} 
A(V(A)) \cap V(A) = \emptyset.
\ee

The last fact holds for every $m \geq 2$. Indeed, if $A = P/Q$ belongs to $\operatorname{Rat}_m \setminus N$, then the critical values of $A$ are finite and distinct, and are given by the zeros of the polynomial $R(t)$ defined by equality \eqref{res}. Consequently, the zero set of the polynomial $S(t)$ given by
\be \l{l}
S(t) = \operatorname{Res}_{2m-2, m, z} \bigl( R(z),\, P(z) - Q(z)t \bigr)
\ee
coincides with the set of finite values attained by $A$ at its critical values.
Therefore, if $Z \subset \mathbb{CP}^{2m+1}$ denotes the projective hypersurface defined by
\be \l{z}
Z : \operatorname{Res}_{2m-2, 2m-2, t} \bigl( R(t), S(t) \bigr) = 0,
\ee
then for every rational function $A \in \operatorname{Rat}_m \setminus (N \cup Z)$, condition \eqref{av} holds. 
The non-emptiness of this Zariski open subset is readily verified. 

Finally, for $m=2$, the theorem follows from Lemma \ref{d2}. 
Indeed, it is well known that a general rational function of degree $m \geq 2$ is not an ordinary Latt\`es map. On the other hand, if $A$ is a generalized Latt\`es map of degree  two that is not an ordinary Latt\`es map, then Lemma \ref{d2} implies that either \eqref{stt} holds, or 
$A: \f O \rightarrow \f O$ is a minimal holomorphic map for some orbifold $\f O$ with $\nu(\f O) = \{2,2\}$. In the latter case, it is easy to see that if $z_1$ and $z_2$ are points with $\nu(z_1)=\nu(z_2)=2$, then, up to relabeling,  $A(z_1)=A(z_2)=z_1$, where $z_1$ is not a critical value of $A$, while $z_2$ is a critical value whose preimage is not a ramified point of $\f O$.
Thus, in this case 
\be 
A(V(A)) \cap \operatorname{Fix}(A) \neq \emptyset,
\ee
and it is easy to see that for general rational functions the last condition, similar to \eqref{stt}, does not hold.
 \qed 

\subsection{The triviality of $\operatorname{Aut}_\infty(A)$ and $G(A)$ for general $A$}
\label{subsec:generic-iterate-automorphisms}
We start by filling a gap in the following result from \cite{ite} (Lemma 3.10). In \cite{ite}, it was erroneously formulated for $m\ge3$,
while it is true (and used) for $m\ge4$. Moreover, 
the proof in \cite{ite} has a gap: after clearing denominators in 
the identity $\alpha\circ A\circ\beta=A$, it considers a projective 
incidence variety 
\[ 
\mathcal I_m\subset 
\mathbb P^3\times\mathbb P^{2m+1}\times\mathbb P^3. 
\] 
However, $\mathcal I_m$ contains the trivial section 
$\{\mathrm{id}\}\times\mathbb P^{2m+1}\times\{\mathrm{id}\}$, so its 
projection to the middle factor is the whole space. Thus, the closedness 
of the projection  does not yield a proper closed subset 
containing the functions with non-trivial $G(A)$. 


For a finite subset $S\subset\mathbb P^1$, write
\[
\operatorname{Stab}(S)
=\{\mu\in\operatorname{PGL}_2(\mathbb C):\mu(S)=S\}.
\]

\begin{lemma}\label{again}
For a general rational function $F$ of degree $m\ge4$,
the group $G(F)$ is trivial.
\end{lemma}

\begin{proof}
Set $r=2m-2$,  and let
$\mathcal C_r\subset\operatorname{Sym}^r(\mathbb P^1)$
be the open subset consisting of reduced divisors,
identified with unordered configurations of $r$ distinct points. Note that $m\geq 4$ implies $r\ge 6$.   

We first show that the configurations for which no two
distinct four-point subsets are M\"obius equivalent form
a non-empty Zariski open subset $\Omega_r\subset\mathcal C_r$.
Consider ordered configurations $(z_1,\ldots,z_r)$.
For two distinct four-element index sets
$\{i_1,i_2,i_3,i_4\}$ and $\{j_1,j_2,j_3,j_4\}$,
and a permutation $\pi\in S_4$, the equality 
of cross-ratios 
\[
(z_{i_1},z_{i_2};z_{i_3},z_{i_4})
=
(z_{j_{\pi(1)}},z_{j_{\pi(2)}};
 z_{j_{\pi(3)}},z_{j_{\pi(4)}})
\]
defines a proper closed subset of the space of ordered
configurations, since the two
sides depend on different sets of variables.
The union of these finitely many closed subsets is also proper and is invariant under $S_r$, so its image under
the finite quotient map to $\mathcal C_r$ is a proper
closed subset. 

Let $C\in\Omega_r$ and let $\mu$ be a M\"obius transformation
such that $\mu(C)=C$. By the definition of $\Omega_r$,
we have $\mu(S)=S$ for every four-point subset $S\subset C$.
Since $r\ge6$, this implies that $\mu=\mathrm{id}$.
Indeed, if $\mu(x)\ne x$ for some $x\in C$, there exists
a four-point subset $S\subset C$ containing $x$ but not
$\mu(x)$, contradicting $\mu(S)=S$.

By Theorem~1.2 of~\cite{Goldberg}, the critical-divisor morphism
\[
\operatorname{Crit}_m:\operatorname{Rat}_m
\longrightarrow\operatorname{Sym}^r(\mathbb P^1),
\qquad
F\longmapsto\sum_{z\in\mathbb P^1}(\deg_zF-1)[z],
\]
is dominant. Therefore,
$\operatorname{Crit}_m^{-1}(\Omega_r)$ is a non-empty
Zariski open subset of $\operatorname{Rat}_m$.
For every $F$ in this subset, the critical divisor is
reduced and $\operatorname{Stab}(C(F))=\{\mathrm{id}\}$.
Since  $\sigma\in G(F)$ implies 
$\sigma(C(F))=C(F)$, we conclude that $G(F)=\{\mathrm{id}\}.$
\end{proof}

For a rational function $A$, set 
\[
K(A)=\operatorname{Stab}(C(A))\cap\operatorname{Stab}(V(A)).
\]

\begin{lemma}\label{lama}
For every $m\geq3$, there exists a non-empty Zariski open
subset \linebreak
$U_m\subset\operatorname{Rat}_m$  such that $K(A)=\{\mathrm{id}\}$ for every $A\in U_m$.
\end{lemma}

\begin{proof}
Put $r=2m-2$, and keep the notation $W,R,L,D$
introduced above. Let $S_m$ be the set of all simple rational functions
$A\in\operatorname{Rat}_m$ such that
$A\notin L\cup D$ and
\[
\operatorname{Res}_{r,r,t}(W(t),R(t))\ne0.
\]
Then $S_m$ is a Zariski open subset of
$\operatorname{Rat}_m$ and for $A\in S_m$,   by construction,  $C(A)$ and $V(A)$ are disjoint sets, each consisting
of $r$ distinct complex numbers. 

For each ordered triple $T=(a,b,c)$ of distinct points
of $C(A)$, let
\[
f_T(z)=\frac{(b-c)(z-a)}{(b-a)(z-c)} 
\]
be the unique M\"obius transformation sending
$a,b,c$ to $0,1,\infty$, respectively. Consider the polynomial
\[
H_A(t)=\prod_{T=(a,b,c)}\prod_{v\in V(A)}
\bigl((b-a)(v-c)t-(b-c)(v-a)\bigr),
\]
where the first product is over all ordered triples
of distinct points of $C(A)$. Since $b-a\ne0$ and $v-c\ne0$,
each factor of $H_A(t)$ has degree one and its root is $f_T(v)$.

We first show that the condition that $H_A$ has distinct roots, or equivalently that 
distinct pairs $(T,v)$ give distinct numbers $f_T(v)$, implies
$K(A)=\{\mathrm{id}\}$. Suppose that $\mu\in K(A)$.
By definition, $\mu$ preserves both $C(A)$ and $V(A)$.
Hence, for any $T=(a,b,c)$ and $v\in V(A)$, the triple
\[
\mu(T)=(\mu(a),\mu(b),\mu(c))
\]
also consists of distinct critical points, and $\mu(v)$
is a critical value. Thus both pairs $(T,v)$ and
$(\mu(T),\mu(v))$ occur in the product defining $H_A$.

The transformations $f_T$ and $f_{\mu(T)}\circ\mu$
both send $a,b,c$ to $0,1,\infty$, respectively.
They therefore coincide, giving
\be \l{mr} 
f_{\mu(T)}(\mu(v))=f_T(v).
\ee 
If $\mu\ne\mathrm{id}$, then $\mu(T)\ne T$, since otherwise
$\mu$ would fix three distinct points, 
and hence $H_A$ has a multiple root by \eqref{mr}.  
Thus, if $H_A$ has distinct roots, $\mu=\mathrm{id}$.

The coefficients of $H_A$ are symmetric polynomials separately
in the roots of $W$ and in the roots of $R$, and hence can be
expressed in terms of their coefficients. Therefore, the condition
$\operatorname{Discr}(H_A)=0$ is algebraic in the coefficients
$a_i,b_i$ of the polynomials $P$ and $Q$ in~\eqref{coe}. Moreover, 
since simultaneous rescaling of all $a_i,b_i$ leaves
$H_A$ unchanged, the set 
\[
U_m=\{A\in S_m:\operatorname{Discr}(H_A) \ne0\}
\]
is a Zariski open subset of $\operatorname{Rat}_m$,
and $K(A)=\{\mathrm{id}\}$ for every $A\in U_m$.

To prove non-emptiness of $U_m$, fix a simple rational
function $A_0\notin L\cup D$ and choose
$p\in\mathbb C\setminus C(A_0)$ such that the numbers
$f_T(p)$ are pairwise distinct for different ordered
triples $T$ of distinct points of $C(A_0)$.
Consider the functions 
\[
A_\varepsilon=p+\varepsilon A_0,\qquad \varepsilon\ne0.
\] 
The critical points of $A_{\varepsilon}$ 
coincide with those of $A_0$, and its critical values
are $p+\varepsilon v$, where $v\in V(A_0)$.
Thus, each $A_{\varepsilon}$ is simple and $A_\varepsilon\notin L\cup D$. Moreover,   
for sufficiently small nonzero $\varepsilon$, the critical
values avoid $C(A_0)$, and hence $A_\varepsilon\in S_m$.

Since
\[
f_T(p+\varepsilon v)\longrightarrow f_T(p)
\quad\text{as} \quad \varepsilon\to0,
\]
values corresponding to different triples $T$ remain
distinct for sufficiently small $\varepsilon$.
For a fixed $T$, distinct critical values $v$ give
distinct values $f_T(p+\varepsilon v)$ by the
injectivity of $f_T$. Thus $H_{A_\varepsilon}$ has
distinct roots for sufficiently small nonzero
$\varepsilon$, and hence $A_\varepsilon\in U_m$.
This proves the lemma.
\end{proof}

\begin{corollary}\label{cor:generic-iterate-automorphisms}
For every $m\ge2$, there exists a non-empty Zariski open
subset $W_m\subset\operatorname{Rat}_m$ such that
$\operatorname{Aut}_\infty(A)=\{\mathrm{id}\}$
for every $A\in W_m$.
\end{corollary}

\begin{proof}
For $m=2$, the assertion is Corollary~\ref{c52}.
Suppose that $m\ge3$. Let $W_m$ be the intersection
of the open subset $U_m$ from Lemma~\ref{lama} with
a non-empty Zariski open subset on which the conclusion
of Theorem~\ref{t1} holds. 

Fix $A\in W_m$. If $\mu\in\operatorname{Aut}(A)$,
then $A\circ\mu=\mu\circ A$ implies that $\mu$
preserves both $C(A)$ and $V(A)$.
Thus $\mu\in K(A)=\{\mathrm{id}\}$.
On the other hand, if
$\mu\in\operatorname{Aut}(A^{\circ n})$
for some $n\ge2$, then applying Theorem~\ref{t1}
to the decomposition
\[
A^{\circ n}
=(\mu^{-1}\circ A)\circ
\underbrace{A\circ\cdots\circ A}_{n-2\ \text{factors}}
\circ(A\circ\mu),
\]
we obtain M\"obius transformations $\nu,\delta$ such that
\[
\mu^{-1}\circ A=A\circ\nu,\qquad
A\circ\mu=\delta\circ A.
\]
Hence, again, $\mu\in K(A)=\{\mathrm{id}\}$.
\end{proof} 

\subsection{The equation $A^{\circ n}=B^{\circ n}$ for general $A$}

For $m\geq 4$, the following result follows from Theorem~1.3
in~\cite{ite}. Below, we derive it for all $m\geq 2$ from Corollary~\ref{cor1} and
Corollary~\ref{cor:generic-iterate-automorphisms}.

\begin{theorem}  \l{it}
For a general rational function $A$ of degree $m\geq 2$, the
equality
\begin{equation} \l{ic1} 
A^{\circ n}=B^{\circ n}
\end{equation}
for some rational function $B$ of degree $m$ and integer
$n\geq 1$ implies that $B=A$.
\end{theorem}

\begin{proof}
We show that the conclusion of the theorem is true for every $A$ in the intersection of non-empty Zariski open subsets
on which the conclusions of Corollary~\ref{cor1} and
Corollary~\ref{cor:generic-iterate-automorphisms} hold.  The case $n=1$ is immediate, so assume that
$n\geq 2$.

Applying Corollary~\ref{cor1} to the decomposition
$$
A^{\circ n}=B^{\circ(n-1)}\circ B
$$
and using $\deg B=\deg A=m$, we obtain a M\"obius
transformation $\delta$ such that
\begin{equation}
\label{eq:common-iterate-end-factors}
B^{\circ(n-1)}=A^{\circ(n-1)}\circ\delta^{-1}, \qquad 
B=\delta\circ A.
\end{equation}
Consequently,
\[
A^{\circ n}
=B\circ B^{\circ(n-1)}
=\delta\circ A^{\circ n}\circ\delta^{-1}.
\]
Thus $\delta\in\operatorname{Aut}(A^{\circ n})$, and
Corollary~\ref{cor:generic-iterate-automorphisms} gives
$\delta=\mathrm{id}$. Hence $B=A$.
\end{proof}

\subsection{Invariant curves} 
Our proof of Theorem~\ref{cur} relies on the following two results on semiconjugate rational functions and invariant curves for endomorphisms $(A_1, A_2)$ of $(\mathbb{C}\mathbb{P}^1)^2$, in the case where the functions involved are not generalized Latt\`es maps.

The first result is a corollary of the classification of semiconjugate rational functions (see \cite{ic}, Proposition 3.3):

\bt \l{tt}
Let $A$ and $B$ be rational functions of degree at least two, and let $X$ be a rational function of degree at least one such that
$
A \circ X = X \circ B.
$
Assume that $A$ is not a generalized Latt\`es map. Then there exist a rational function $Y$ and an integer $d \geq 0$ such that
$
X \circ Y = A^{\circ d}.
$ \qed 
\et

The second result is a corollary of the classification of  invariant curves (see \cite{ic}, Theorem 1.1):

\bt \l{1}
Let $A_1$ and $A_2$ be rational functions of degree at least two that are not generalized Latt\`es maps, and let $C$ be an irreducible algebraic curve in $(\C\P^1)^2$ that is invariant under $(A_1, A_2)$ and is not a vertical or horizontal line. Then there exist rational functions $X_1$, $X_2$, $Y_1$, $Y_2$, and $B$ such that:
\begin{enumerate}
\item The diagram
\be \l{xx}
\begin{CD}
(\C\P^1)^2 @>(B,B)>> (\C\P^1)^2 \\
@V (X_1, X_2) VV @VV (X_1, X_2) V \\
(\C\P^1)^2 @>(A_1, A_2)>> (\C\P^1)^2
\end{CD}
\ee
commutes.
\item The equalities
\[
X_1 \circ Y_1 = A_1^{\circ d}, \quad X_2 \circ Y_2 = A_2^{\circ d}
\]
hold for some integer $d \geq 0$.
\item The map $t \mapsto (X_1(t), X_2(t))$ parametrizes the curve $C$.
\qed 
\end{enumerate}
\et

Theorems~\ref{tt} and~\ref{lats} combined with results of \cite{ite} imply the following statement.

\bt \label{t4}  
For a general rational function \( A \) of degree \( m \geq 2 \), the following holds:  
whenever \( B, X \) are non-constant rational functions such that the diagram  
\be \label{ii2}
\begin{CD}
\C\P^1 @>B>> \C\P^1 \\
@V X VV @VV X V\\ 
\C\P^1 @>A^{\circ r}>> \C\P^1
\end{CD}
\ee
commutes for some integer \( r \geq 1 \), there exist a Möbius transformation \( \mu \) and an integer \( l \geq 0 \) such that
\be \l{bo} 
X = A^{\circ l} \circ \mu, \qquad B = \mu^{-1} \circ A^{\circ r} \circ \mu.
\ee 
\et

\pr Since the statement of the theorem holds for simple rational functions
of degree $m\geq4$ by Theorem~1.4 of~\cite{ite}, it remains to consider
the cases $m=2$ and $m=3$.
Let $U\subset\operatorname{Rat}_m$ be a non-empty Zariski open subset
such that every $A\in U$ is simple and satisfies the conclusions of
Theorem~\ref{lats}, Corollary~\ref{cor1}, and, in the case $m=2$,
Corollary~\ref{c41}.

Fix $A\in U$ and suppose that diagram~\eqref{ii2} commutes.
If $\deg X=1$, take $\mu=X$ and $l=0$. Both equalities
in~\eqref{bo} then follow directly from~\eqref{ii2}.
Henceforth assume that $\deg X\geq2$. Comparing degrees
in~\eqref{ii2}, we obtain $\deg B=m^r\geq2$.

Since a rational function is a generalized Latt\`es map if and only
if some iterate of it is (see \cite{ic}, Section~2.3), and $A\in U$
is not a generalized Latt\`es map, it follows that $A^{\circ r}$ is
also not a generalized Latt\`es map. Hence, by Theorem~\ref{tt},
there exist a rational function $Y$ and an integer $d\geq0$ such that
\[
X\circ Y=A^{\circ rd}.
\]
Since $\deg X\geq2$, necessarily $d\geq1$.
If $\deg Y=1$, then $X=A^{\circ rd}\circ Y^{-1}$.
Otherwise, Corollary~\ref{cor1} applies to this decomposition.
Thus, in either case, there exist a M\"obius transformation $\mu$
and an integer $l\geq1$ such that
\be \l{ob}
X=A^{\circ l}\circ\mu.
\ee
This is the first equality in~\eqref{bo}. Substituting~\eqref{ob}
into diagram~\eqref{ii2}, we obtain
\[
A^{\circ r}\circ A^{\circ l}\circ\mu
=A^{\circ l}\circ\mu\circ B,
\]
and hence
\be \l{gir}
A^{\circ l}\circ A^{\circ r}
=A^{\circ l}\circ\mu\circ B\circ\mu^{-1}.
\ee

In the case $m=3$, applying Theorem~\ref{goo1} repeatedly
to~\eqref{gir}, we obtain
\[
A^{\circ r}=\mu\circ B\circ\mu^{-1},
\]
which gives the second equality in~\eqref{bo}.

In the case $m=2$, applying Corollary~\ref{cor1} to~\eqref{gir}
and comparing degrees, we obtain a M\"obius transformation $\delta$
such that
\[
A^{\circ l}=A^{\circ l}\circ\delta,
\qquad
\mu\circ B\circ\mu^{-1}=\delta^{-1}\circ A^{\circ r}.
\]
The first equality implies that $\delta\in\Sigma_\infty(A)$,
and Corollary~\ref{c41} gives $\delta\in\Sigma(A)$.
Hence, for $\mu'=\delta\circ\mu$, we have 
\[
X=A^{\circ l}\circ\mu'
\]
and 
\[
B=\mu^{-1}\circ\delta^{-1}\circ A^{\circ r}\circ\mu
=\mu'^{-1}\circ A^{\circ r}\circ\mu'.
\]
Thus both equalities in~\eqref{bo} hold with $\mu'$ in place of $\mu$,
and the theorem also follows for $m=2$. \qed

\vskip 0.2cm
\noindent{\it Proof of Theorem \ref{cur}.} By the results of \cite{ite}, it suffices to prove the theorem for \( m = 2 \) or \( 3 \). So, let \( m \) be either of these values, and let \( U \subset \operatorname{Rat}_m \) be a non-empty Zariski  open subset for which the conclusions of Theorem \ref{lats},   Theorem \ref{it},   and Theorem \ref{t4} hold.

Suppose \( A_1, A_2 \in U \) and 
\be \l{lass} 
(A_1, A_2)^{\circ d}(C) = C, \quad d \geq 1.
\ee
Then Theorem \ref{1} and Theorem \ref{t4} imply that \( C \) is parametrized by
\[
t \mapsto \bigl( (A_1^{\circ d_1} \circ \beta)(t), \, (A_2^{\circ d_2} \circ \alpha)(t) \bigr)
\]
for some integers \( d_1, d_2 \geq 0 \) and M\"obius transformations \( \alpha, \beta \) such that
\[
\beta^{-1} \circ A_1^{\circ d} \circ \beta = \alpha^{-1} \circ A_2^{\circ d} \circ \alpha.
\]
Moreover, without loss of generality, we may assume that \( \beta \) is the identity map, which implies
\[
A_1^{\circ d} = \alpha^{-1} \circ A_2^{\circ d} \circ \alpha = (\alpha^{-1} \circ A_2 \circ \alpha)^{\circ d}.
\]
By Theorem \ref{it}, this yields
\be \l{kosh1}
A_2 = \alpha \circ A_1 \circ \alpha^{-1}.
\ee
 Thus, the parametrization above becomes
\[
t \mapsto \bigl( A_1^{\circ d_1}(t), \, \alpha \circ A_1^{\circ d_2}(t) \bigr).
\]

If \( d_1 \leq d_2 \), then this parametrization reduces to
\be \l{par1}
t \mapsto \bigl(t,\, (\alpha \circ A_1^{\circ s})(t)\bigr),
\ee where $s=d_2-d_1.$
On the other hand, if \( d_1 > d_2 \), then \( C \) is parametrized by
\[
t \mapsto \bigl( A_1^{\circ s}(t), \, \alpha(t) \bigr),
\] where $s=d_1-d_2$, 
and hence also by
\be \l{par2}
t \mapsto \bigl((A_1^{\circ s} \circ \alpha^{-1})(t),\, t\bigr).
\ee

To finish the proof, it remains to observe that if \eqref{kosh1} holds, then a curve $C$ parametrized by \eqref{par1} or \eqref{par2} 
is $(A_1,A_2)$-invariant, which can be verified by a direct calculation. 
Indeed, in the first case, if the point $(x_0, y_0)$ belongs to $C$, then 
\[
y_0 = (\alpha \circ A_1^{\circ s})(x_0),
\] 
which implies that 
\[
A_2(y_0) = A_2\bigl((\alpha \circ A_1^{\circ s})(x_0)\bigr) = (\alpha \circ A_1^{\circ s})(A_1(x_0)).
\] 
Thus, the point $(A_1(x_0), A_2(y_0))$ also belongs to $C$. In the second case, the proof is similar. \qed

\section*{acknowledgements}
The author is grateful to Geng-Rui Zhang for his comments on the paper.


\begin{thebibliography}{99}

\bibitem{mm}
L.~DeMarco, N.~M.~Mavraki,
\textit{Multiplier spectrum of maps on $\P^1$:
theorems of Ji--Xie and Ji--Xie--Zhang},
notes prepared for lectures at Harvard and Toronto
in November 2023, revised April 2024.

\bibitem{f2}
M.~Fried,
\textit{The field of definition of function fields and a problem
in the reducibility of polynomials in two variables},
Illinois J. Math. 17 (1973), no.~1, 128--146.

\bibitem{f3}
M.~Fried,
\textit{Arithmetical properties of function fields (II).
The generalized Schur problem},
Acta Arith. 25 (1974), no.~3, 225--258.

\bibitem{Goldberg}
L.~R.~Goldberg,
\textit{Catalan numbers and branched coverings by the Riemann sphere},
Adv. Math. 85 (1991), no.~2, 129--144.

\bibitem{go}
R.~Gotou,
\textit{Dynamical systems of correspondences on the projective
line II: degrees of multiplier maps},
arXiv:2309.15404 (2023).

\bibitem{har}
J.~Harris,
\textit{Galois groups of enumerative problems},
Duke Math. J. 46 (1979), no.~4, 685--724.

\bibitem{hug}
V.~Huguin,
\textit{Moduli spaces of polynomial maps and multipliers
at small cycles},
arXiv:2412.19335 (2024).

\bibitem{hut}
B.~Hutz, M.~Tepper,
\textit{Multiplier spectra and the moduli space of degree 3
morphisms on $\P^1$},
JP J. Algebra Number Theory Appl. 29 (2013), no.~2, 189--206.

\bibitem{xie}
Z.~Ji, J.~Xie,
\textit{The multiplier spectrum morphism is generically injective},
J. Eur. Math. Soc. (2025), published online first.

\bibitem{ms}
A.~Medvedev, T.~Scanlon,
\textit{Invariant varieties for polynomial dynamical systems},
Ann. of Math. (2) 179 (2014), no.~1, 81--177.

\bibitem{mil1}
J.~Milnor,
\textit{Geometry and dynamics of quadratic rational maps},
Experiment. Math. 2 (1993), no.~1, 37--83.

\bibitem{mil2}
J.~Milnor,
\textit{On Latt\`es maps},
in \textit{Dynamics on the Riemann Sphere:
A Bodil Branner Festschrift},
P.~G.~Hjorth and C.~L.~Petersen (eds.),
European Mathematical Society, 2006, pp.~9--43.

\bibitem{pak}
F.~Pakovich, \textit{Functional equations in branched covering maps}, arXiv:2609.23613.  



\bibitem{pq}
F.~Pakovich,
\textit{Algebraic curves $P(x)-Q(y)=0$ and functional equations},
Complex Var. Elliptic Equ. 56 (2011), no.~1--4, 199--213.

\bibitem{semi}
F.~Pakovich,
\textit{On semiconjugate rational functions},
Geom. Funct. Anal. 26 (2016), 1217--1243.

\bibitem{lattes}
F.~Pakovich,
\textit{On generalized Latt\`es maps},
J. Anal. Math. 142 (2020), no.~1, 1--39.

\bibitem{ic}
F.~Pakovich,
\textit{Invariant curves for endomorphisms of $\P^1\times\P^1$},
Math. Ann. 385 (2023), no.~1--2, 259--307.

\bibitem{low}
F.~Pakovich,
\textit{Lower bounds for genera of fiber products},
Israel J. Math. 269 (2025), 475--500.

\bibitem{sym}
F.~Pakovich,
\textit{On symmetries of iterates of rational functions},
Ann. Sc. Norm. Super. Pisa Cl. Sci. (5) 26 (2025),
no.~3, 1677--1701.

\bibitem{ite}
F.~Pakovich,
\textit{On iterates of rational functions with maximal
number of critical values},
J. Anal. Math. 156 (2025), no.~1, 213--251.

\bibitem{we}
L.~W.~West,
\textit{The moduli space of cubic rational maps},
arXiv:1408.3247 (2014).

\bibitem{ye}
H.~Ye,
\textit{Rational functions with identical measure of maximal entropy},
Adv. Math. 268 (2015), 373--395.

\end{thebibliography}
\end{document}